# A New Continuum-Based Thick Shell Finite Element for Soft Biological Tissues in Dynamics: Part 2 – Anisotropic Hyperelasticity and Incompressibility Aspects

Bahareh Momenan, Ph.D.
*Department of Mechanical Engineering, University of Ottawa,*
*161 Louis Pasteur, Ottawa, K1N 6N5 Ontario, Canada*
baharehmomenan@gmail.com

Michel R. Labrosse, Ph.D.
*Department of Mechanical Engineering, University of Ottawa,*
*161 Louis Pasteur,* Colonel By Hall A-213, *Ottawa, K1N 6N5 Ontario, Canada*
labrosse@eng.uottawa.ca

In a companion article (Part 1), we presented the development of a thick continuum-based (CB) shell finite element (FE) based on Mindlin-Reissner theory. We verified the accuracy, efficiency and locking insensitivity of the element in modeling large 3D deformations, using linear elastic material properties. In the present article, we developed and implemented the kinetics description, within the updated Lagrangian (UL) formulation, of anisotropic incompressible hyperelastic constitutive relations that enable the CB shell FE to accurately model very large 3D strains and deformations. Specifically, we developed the measures of deformation in the lamina coordinate system, presented three techniques to model nonlinear hyperelastic strains, and enabled the direct enforcement of incompressibility and of the zero normal stress condition without using a penalty factor or a Lagrange multiplier.

Moving towards the application of the present work to the biomedical realm, we performed multiple experiments concerning mechanical behavior of rubber-like materials and soft biological tissues in different geometries and loading conditions. Excellent agreements between the present FE results and the analytical and/or experimental data proved the CB shell FE combined with the present constitutive techniques to be a highly reliable and efficient tool for modeling, analyzing, and predicting mechanical behavior of soft biological tissues.

*Keywords:* Anisotropic hyperelastic incompressible; soft tissue dynamics; updated Lagrangian formulations.

## 1. Introduction

It has been shown that the development of many diseases (such as atherosclerosis, heart failure, and asthma) may be associated with changes in cell mechanics, extracellular matrix structure, or mechanotransduction (i.e., the mechanisms by which cells sense and respond to mechanical signals) [1]. Therefore, it is important to be able to properly understand and simulate the mechanical behaviour of biological tissues. Furthermore, accurate tissue deformation information and force feedback are also needed for the simulation of surgical procedures.

While hard tissues (e.g. bone) benefit from similarities with classical engineering materials, rubber-like materials and soft biological tissues (e.g. blood vessels, aortic valve leaflets) require the use of advanced tools, such as finite element (FE) methods and strain energy functions (hyperelasticity), for their modeling and simulation. Assuming that a strain energy function and material constants are properly determined, commercial FE codes mostly provide users only with volume finite elements such as bricks [2–5]. Due to the limitations of brick elements and some existing shell elements to model large 3D deformations and large 3D strains of biological soft tissues, as detailed in a companion article (Part 1), we developed a new 9-noded quadrilateral thick continuum-based (CB) shell FE to handle such problems, and presented the kinematic derivations and the constitutive formulation of linear elasticity. However, to make full use of this element for the simulation of rubber-like materials and soft biological tissues, the appropriate constitutive models must be developed.

In the literature, hyperelastic constitutive relations are developed by differentiating the strain energy function $W$ with respect to strains within the total Lagrangian (TL) formulation, in which the initial (undeformed) configurations is considered as the reference to measure loads and deformations [6,7]. However, when large 3D deformations and strains are considered, the updated Lagrangian (UL) formulation has been argued to be more efficient and more accurate [8]; it is used herein.

In the case of membrane or shell elements, the incompressibility condition is typically added by modifying the strain energy function such that: $\widetilde{W} = W(\bar{\bar{F}}) - p(\det\bar{\bar{F}} - 1)$. In this formulation, deformation gradient $\bar{\bar{F}}$ is obtained due to the distortional effects only, and Lagrange multiplier $p$ enforcing $\det\bar{\bar{F}} = 1$ is an additional unknown [6,7,9,10]. Authors of reference [10] solve for the Lagrange multiplier within the context of thin shells (i.e.

Kirchhoff-Love based shell formulations, $S_{33} = S_{13} = S_{23} = 0$), where the transverse shears are neglected, and the transverse normal strain $E_{33}$, which cannot be neglected in the case of large strains, is statically condensed using the plane stress condition. This constitutes a 2D specialization within the framework of a 3D material model. However, according to reference [11], in a 2D material model, a significant part of the 3D constitutive law is missing, because $W$ does not depend on $E_{33}$, and thus $\partial W/\partial E_{33} = 0$. Therefore, a 2D material model does not describe the three-dimensional properties, and does not distinguish between compressible and incompressible materials.

From a different point of view, the additive decomposition of the strain-energy function into volumetric and deviatoric (distortional) parts that has been embedded into most 3D FE implementations to date, is essentially an isotropic condition, and is appropriate for pure hydrostatic tension only [12]. Thus, its generalization to anisotropic materials is not adequate, considering the corresponding physics. The arguments above suggest that, despite its intuitive appeal, it should not be employed when modeling nonlinear, anisotropic materials which are characterized by infinitesimal volume changes when deformed. Certainly, its equivalent formulation in terms of stresses (application of the zero normal stress condition) does not seem a natural or appropriate constitutive assumption to make when modeling nearly incompressible or fully incompressible anisotropic materials. Thus, a formulation of the theory that accounts for infinitesimal volume changes in a physically realistic way is badly needed.

The objective in the present work is to develop the kinetics description and constitutive relations for anisotropic near-incompressible hyperelastic materials typical of rubber-like materials and soft biological tissue dynamics. When combined with the CB shell FE presented in Part 1, and the updated Lagrangian (UL) formulation, the goal is to allow for very large 3D deformations and very large 3D strains using a coarse mesh.

To do so, we developed the measures of deformation in the lamina coordinate system (Section 2.3) as opposed to the global ones typically used in literature, and further evaluated the constitutive relations in the lamina coordinate system. We developed a new technique to implement hyperelastic nonlinear anisotropic incompressible constitutive relations in the lamina coordinate system, within the total UL formulation (Technique 1, Section 2.4.1). In this technique, we usedd a direct method to enforce the zero normal stress condition through incompressibility, without decomposing the strain energy function into volumetric and distortional components. In addition, we developed two other techniques that enable very large strains using a constant elasticity tensor such as found in linear elasticity. The latter techniques, especially applicable for rubber-like materials, are computationally more efficient than Technique 1, and circumvent the need of evaluating material constants for strain energy functions. Technique 2 is based on the incremental decomposition of stresses (Section 2.4.2) within the total UL formulation, and Technique 3 is based on the linearization of the constitutive relation (Section 2.4.3) within the incremental UL formulation.

To verify the accuracy and efficiency of the present CB shell FE in very large 3D strains and very large 3D deformations, and the robustness of the three techniques to evaluate the hyperplastic constitutive relations in the lamina coordinate system within the incremental and total UL formulations, we performed multiple experiments concerning rubber-like materials and soft biological tissues with different geometries and loading conditions.

## 2.  Methods

### 2.1.  *Geometric and kinematic description*

The reader is referred to Part 1 (Section 2.1) for the geometric and kinematic description of the present 9-noded CB shell FE.

In the UL formulation, the equilibrium is considered at the current configuration, and loads and displacements are either measured from the initial/undeformed configuration (total UL) [6,13,14], or from the previously deformed configuration (incremental UL) [8]. Herein, the left superscript $\tau$ denotes the current configuration, and the left subscript $\beta$ denoting the reference configuration, is replaced by 0 in the total UL formulation, and by $\tau - \Delta\tau$ in the incremental UL formulation.

### 2.2.  *Coordinate systems*

To enable large in-plane strains and large bending deformations/rotations, we introduced two independent coordinate systems: one that remains normal to the surface of the shell (lamina coordinate system: $\vec{e}^l$); and a second that is corotational with the fiber (fiber coordinate system: $\vec{e}^f$) [13,15]. As discussed in Part 1, the independency of the fiber from the normal of the shell allows for large bending deformations and rotational strains without creating artificial stiffening, and thus prevents shear and membrane locking.

The derivations of the lamina coordinate systems presented in references [14–18] and Part 1 is general and efficient for any type of analysis. However, we discovered that the fiber coordinate system algorithm proposed in reference [15] is limited to deflections (bending deformations) smaller than 90°. Explicitly, if the bending deformation exceeds 90°, this algorithm results in a flip in the orientation of $\vec{e}_1^f$ and $\vec{e}_2^f$ (Fig. 1, right), and ruins the application of the boundary conditions as well as the evaluation of the force and displacement vectors.

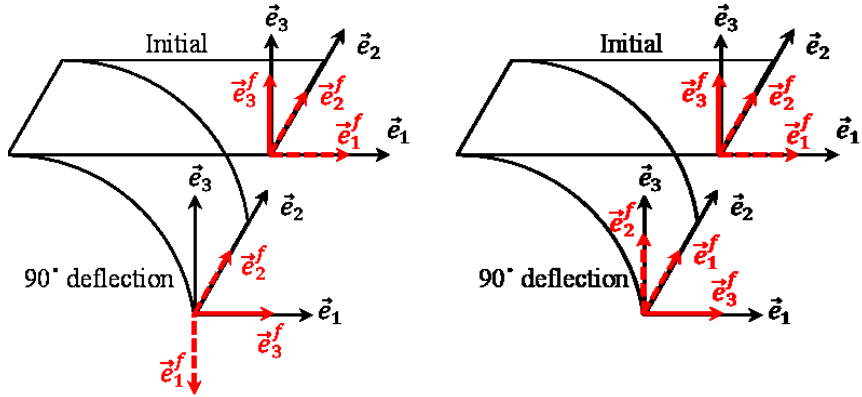

Fig. 1. Left: physically expected rotation of the fiber coordinate system as the element deflects. Right: Fiber coordinate system at the 90° deflection of the element obtained from the algorithm presented in [15], changes orientation from that of the initial configuration.

To preserve the physically expected orientation of the fiber coordinate system in deformations larger than 90 degrees (Fig. 1, left), we propose a more robust and simplified algorithm. Assuming that, in the UL formulation, the increments in rotations ($_{\tau-\Delta\tau}^{\tau}\theta_1^a$ and $_{\tau-\Delta\tau}^{\tau}\theta_2^a$) are smaller than 45 degrees, we modified the algorithm already presented in Appendix C in Part 1 to the following [15]:

*Step 1:* $j = 1$.

Knowing that $^{\tau}\hat{Y}$ is the unit basis vector in the fiber direction in the current configuration, we set:

*Step 2:* $^{\tau}\vec{e}_3^f = {}^{\tau}\hat{Y}$

To prevent the flip discussed above, we used a coordinate system whose rotational displacement, due to deformation, was similar to that of the fiber, namely: the fiber coordinate system of the previous time step ($^{\tau-\Delta\tau}\vec{e}_j^f$), to find:

*Step 3:* $^{\tau}\vec{e}_2^f = ({}^{\tau}\hat{Y} \times {}^{\tau-\Delta\tau}\vec{e}_j^f)/\| {}^{\tau}\hat{Y} \times {}^{\tau-\Delta\tau}\vec{e}_j^f \|$

*Step 4:* $^{\tau}\vec{e}_1^f = ({}^{\tau}\vec{e}_2^f \times {}^{\tau}\hat{Y})$

To initialize the solution scheme (i.e. for the first time step), we set the fiber coordinate system as the already calculated lamina coordinate system.

### 2.3. *Measures of deformation in the lamina coordinate system*

To enable large 3D deformations within one shell element, we developed the measures of deformation in the nodal lamina coordinate system (independent and variable from one node to the other within one element, Part 1) as opposed to the current global ones.

We obtained the deformation gradient of the current configuration with respect to the reference configuration in the lamina coordinate system from (1), where the lamina reference and the lamina current Jacobians ($[{}^{\beta}J^l]$ and $[{}^{\tau}J^l]$, respectively), were defined in Part 1.

$$[{}_{\beta}^{\tau}F^l] = [{}^{\tau}J^l]^T[{}^{\beta}J^l]^{-T} \quad (1)$$

We evaluated the lamina Green-Lagrange strain tensor defined in the current configuration ($\tau$) with respect to the reference configuration ($\beta$), and the lamina Almansi strain tensor in the current configuration from (2) and (3), respectively. In addition, the right and left Cauchy-Green tensors in the lamina coordinate systems are respectively obtained from $[{}_{\beta}^{\tau}C^l] = [{}_{\beta}^{\tau}F^l]^T[{}_{\beta}^{\tau}F^l]$ and $[{}_{\beta}^{\tau}B^l] = [{}_{\beta}^{\tau}F^l][{}_{\beta}^{\tau}F^l]^T$.

$$[{}_{\beta}^{\tau}E^l] = \frac{1}{2}([{}_{\beta}^{\tau}C^l] - [I]) \tag{2}$$

$$[{}_{\tau}\varepsilon^l] = \frac{1}{2}\left([I] - [{}_{\beta}^{\tau}B^l]^{-1}\right) \tag{3}$$

### 2.4. *Lamina constitutive relations*

In case of large strains and large deformations, an important aspect is the implementation of the constitutive relations in the lamina coordinate system. In what follows, we present three techniques for the solution of the large strain problems in the total and incremental UL formulations.

#### 2.4.1. *Technique 1: modeling large strains using anisotropic nonlinear hyperplastic constitutive relations, incompressibility, and application of zero normal stress condition*

In *Technique 1*, we considered the total UL formulation, where loads and displacements are measured from the undeformed configuration [16,19,14], and implemented large strains using hyperelastic material models.

The constitutive relations, namely the material tangent modulus and the second Piola-Kirchhoff stress tensor, are typically obtained by taking the derivatives of the strain energy function with respect to the Green-Lagrange strains. However, to have the constitutive relations in the lamina coordinate system, we employed the lamina Green-Lagrange strains as opposed to the global ones. Thus, the lamina material tangent modulus and the lamina second Piola-Kirchhoff stress tensor are, respectively, obtained from (4) and (5).

$$_0C^l_{ijrs} = \frac{\partial^2 {}_0^{\tau}W}{\partial {}_0^{\tau}E^l_{ij}\, \partial {}_0^{\tau}E^l_{rs}} \tag{4}$$

$$_0^{\tau}S^l_{ij} = \frac{\partial {}_0^{\tau}W}{\partial {}_0^{\tau}E^l_{ij}} \tag{5}$$

As needed for the UL formulation, we transformed the lamina material tangent modulus (4) to the current configuration using the lamina deformation gradients, as opposed to the global ones used in reference [19], such that:

$$_{\tau}C^l_{mnpq} = \frac{{}^{\tau}\rho}{{}^0\rho}\, {}_0^{\tau}F^l_{mi}\, {}_0^{\tau}F^l_{nj}\, {}_0C^l_{ijrs}\, {}_0^{\tau}F^l_{pr}\, {}_0^{\tau}F^l_{qs}. \tag{6}$$

Similarly, we obtained the lamina Cauchy stress tensor by transforming the lamina second Piola-Kirchhoff stress tensor (5) using the lamina deformation gradients:

$$^{\tau}\sigma^l_{sr} = \frac{{}^{\tau}\rho}{{}^0\rho}\, {}_0^{\tau}F^l_{si}\, {}_0^{\tau}S^l_{ij}\, {}_0^{\tau}F^l_{rj}. \tag{7}$$

In these relations, left scripts $\tau$ and 0 denote the current and the initial configurations, and $[{}_0^{\tau}F^l]$ and ${}_0^{\tau}E^l_{ij}$ are obtained by replacing $\beta$ with 0 in (1) and (2), respectively. Note that this technique is effective only when the total UL formulation is used with hyperelasticity, in which the constitutive relations are not linear. Therefore, if this transformation is applied to a material with constant constitutive tensor undergoing large strains (Experiment 2 in Section 3.2), totally different results are obtained [19,14].

We note that, so far, no application of zero normal stress condition and/or incompressibility has been considered in this technique. This is addressed in what follows.

Considering hyperelastic materials, authors of [10] suggested that the zero normal stress condition be applied either into the strain energy function via a Lagrange multiplier enforcing incompressibility, or iteratively for compressible materials. However, the limitations and concerns associated with this approach were discussed earlier. As a remedy, we developed a direct approach for enforcing the zero normal stress condition which works equally well for both incompressible and compressible materials.

Considering the total UL formulation, (2) can be re-written as $[{}_0^\tau C^l] = 2[{}_0^\tau E^l] + [I]$, where $[{}_0^\tau C^l] = [{}_0^\tau F^l]^T[{}_0^\tau F^l]$. Thus, using the properties of the determinant and rearranging gives: $\det(2[{}_0^\tau E^l] + [I]) - (\det[{}_0^\tau F^l])^2 = 0$, which is linear in ${}_0^\tau E^l_{33}$. We used this relation to solve for the lamina Green-Lagrange strain in the normal direction (${}_0^\tau E^l_{33}$). The resulting ${}_0^\tau E^l_{33}$ is a function of the remaining eight components of the Green-Lagrange strain tensor only. Thus, if ${}_0^\tau E^l_{33}$ is substituted in the strain energy function, its derivatives with respect to ${}_0^\tau E^l_{33}$ ((4) and (5)) vanish (i.e. ${}_0^\tau S^l_{33} = {}_0 C^l_{33pq} = {}_0 C^l_{mn33} = 0$), and lamina constitutive relations satisfying the zero normal stress condition are developed. This method is valid for both the incompressible and the compressible materials. The only difference is that we set $\det[{}_0^\tau F^l] = 1$ to enforce incompressibility, or we simply use the calculated value of $\det[{}_0^\tau F^l]$ for compressible materials.

2.4.2. *Technique 2: modeling large strains using constant elasticity tensor in the total UL formulation*

*Technique 2* was designed to enable the use of a constant constitutive tensor (e.g. from linear elasticity) to model large strains in the total UL formulation using the incremental decomposition of the stress components in the lamina coordinate system [14,19]:

$$ {}^{\tau+\Delta\tau}_\tau S^l_{ij} = {}^\tau\sigma^l_{ij} + {}_\tau S^l_{ij}. \tag{8} $$

In this relation, the components of the second Piola-Kirchhoff stress increment tensor in the lamina coordinate system referred to the configuration at time $\tau$ are calculated from:

$$ {}_\tau S^l_{ij} = {}^\tau_\tau C^l_{ijrs}\, {}^{\tau+\Delta\tau}_\tau E^l_{rs} \tag{9} $$

where, ${}^\tau_\tau C^l_{ijrs}$ is the constant constitutive tensor, and the components of the Green-Lagrange strain increment in the lamina coordinate system (${}^{\tau+\Delta\tau}_\tau E^l_{ij}$) are obtained by, respectively, replacing the left super- and subscripts of (2) with $\tau + \Delta\tau$ and $\tau$. Finally, similar to (7), the lamina Cauchy stress at time $\tau + \Delta\tau$ are obtained from: ${}^{\tau+\Delta\tau}\sigma^l_{sr} = \frac{{}^{\tau+\Delta\tau}\rho}{{}^\tau\rho}\, {}^{\tau+\Delta\tau}_\tau F^l_{si}\, {}^{\tau+\Delta\tau}_\tau S^l_{ij}\, {}^{\tau+\Delta\tau}_\tau F^l_{rj}$.

The main advantage of this technique is its relatively simple implementation compared to Technique 1. For instance, assuming that Young's modulus and Poisson's ratio are known for small strain analysis, and that a subroutine to calculate the constitutive relations in small strain analysis has been written, the same may be used for large strain analysis, simply by using the accumulation of the stresses (8) and the current incremental Green-Lagrange strains and the incremental second Piola-Kirchhoff stresses (9). Therefore, evaluation of the appropriate material constants for the hyperelastic strain energy function (as in Technique 1), as well as the transformation of the fourth-order constitutive tensor (6) are avoided. Note that the zero normal stress condition needs to be applied as discussed in Part 1 (Section 2.5).

A similar idea for the accumulation of the stresses and the evaluation of incremental stresses and strains from the small strains material law was used in [19] for hypoelasticity including elastoplasticity, in which ${}^\tau_\tau C^l_{ijrs}$ is defined by the history of Cauchy stresses and the accumulation of the instantaneous plastic strain increments.

2.4.3. *Technique 3 modeling large strains using constant elasticity tensor in the incremental UL formulation*

*Technique 3* consists in the incremental linearization of the constitutive relations to model large strains in a material with a constant constitutive tensor (from linear elasticity). Specifically, we considered the incremental UL formulation, in which displacements and loads are measured from the previous configuration. Thus, only the increment in the stresses (as opposed to the accumulation of the stresses used in Technique 2) were considered. In this technique, we multiplied the constant material tensor (without any transformation) with the linearized strain increments in the lamina coordinate system developed in Part 1 (Section 2.4) to obtain the lamina Cauchy stresses in the current configuration:

$$ \{{}^\tau\sigma^l\} = [{}^\tau_\tau C^l]\{{}_\tau e^l\}. \tag{10} $$

In addition to the advantages mentioned in Technique 2, complete exclusion of any transformations makes this computation the most efficient compared to the previous two techniques.

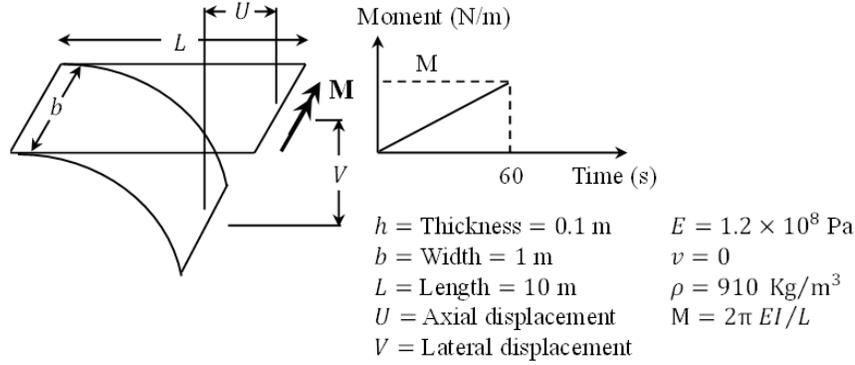

Fig. 2. Geometry, loading condition, and material properties for Experiment 1.

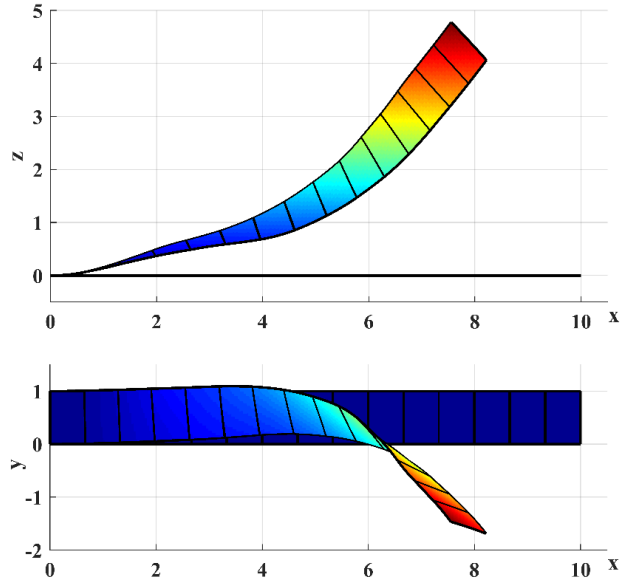

Fig. 3. Failure in the deformation of the cantilever beam subjected to a pure tip bending moment due the limitation of the former fiber coordinate system, as discussed in Section 2.2. Units are in meters.

## 3.  Experiments

To verify the accuracy of the new fiber coordinate system algorithm, and robustness of the three techniques in modeling large 3D deformations and large 3D strains using the incompressibility condition presented above, we performed multiple experiments concerning rubber-like material and soft biological tissues with different geometries and loading conditions. The 9-noded quadrilateral thick CB shell FE presented in Part 1, was employed in these experiments.

### 3.1.  *Experiment 1: Elastic, large pure bending deformations and rotations, small strains*

First, to demonstrate the limitation of the Hughes' fiber coordinate system algorithm, and to evaluate accuracy of the present fiber coordinate system algorithm presented in Section 2.2, we modeled large pure bending of the cantilever beam shown in **Error! Reference source not found.** 2. Failure in large bending deformations due to the limitation of the Hughes' fiber coordinate system is illustrated in Fig. 3. To validate the robustness of the present fiber coordinate system, the FE results are plotted (Fig. 4: moment versus tip axial and transverse displacements) against those of the static analytical solution (reference [20], page 54). Finally, the deformed shape of the cantilever beam at the maximum load configuration obtained from the present fiber coordinate system (Section 2.2) is illustrated in Fig. 5.

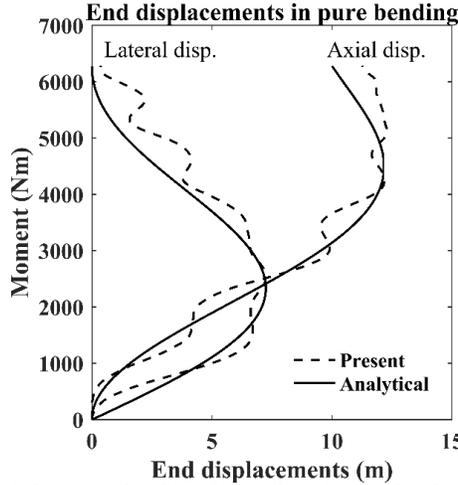

Fig. 4. Moment vs. end displacements in large pure bending of a cantilever beam using the new fiber coordinate system (Section 2.2).

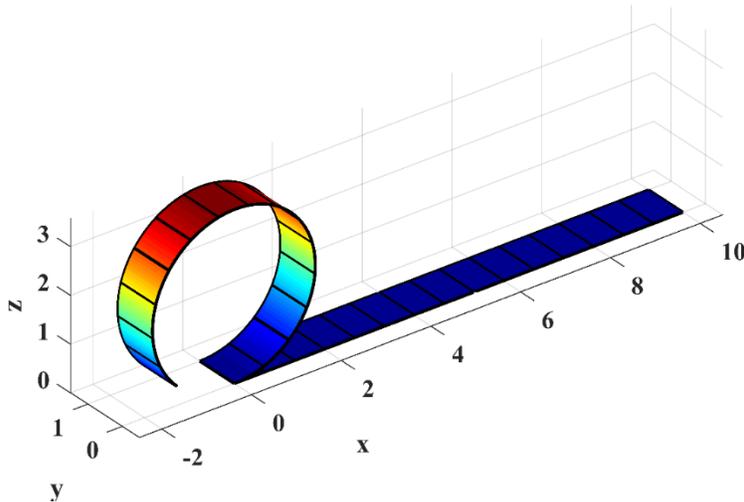

Fig. 5. Deformation of the cantilever beam at the maximum load configuration obtained from the new fiber coordinate system (Section 2.2). Units are in meters.

### 3.2. *Experiment 2: Nonlinear isotropic elastic, large distortions, large strains, verification of Techniques 1-3*

A plane stress large strain analysis of one quarter of a simply supported rubber sheet with a central hole subjected to in-plane edge pressure (Fig. 6, left) was carried out for each of the three techniques presented in Section 2.4. In Technique 1 (Section 2.4.1), we used the hyperelastic incompressible Mooney-Rivlin material model $W = C_1(I_1 - 3) + C_2(I_2 - 3)$, where the material constants $C_1$ and $C_2$, given in [19], were derived from an analytical and experimental investigation done in [21]. For consistency with linear elasticity, relationships can be established between material constants $C_1$ and $C_2$ and Young's modulus and Poisson's ratio. Namely, $E = 2\mu(1 + \nu)$, where $\mu = 2(C_1 + C_2)$, and $\nu$ satisfies the incompressibility assumption. These material constants and properties are provided in Fig. 6, left. Knowing Young's modulus and Poisson's ratio for the material allows for the determination of a constant material elasticity tensor which was used in combination with Technique 2 (Section 2.4.2) and Technique 3 (Section 2.4.3). In these analyses, we employed only two of the present CB shell FE to model one quarter of the rubber sheet (Fig. 6, top-right) and compared the results with those presented in [19] using thirty 4-noded elements per quarter (Fig. 6, top-middle). To have a common ground with the static response analysis presented in [19], we ramp-loaded the plate at a random 400 increments to final load, and presented the load-deformation responses, obtained from the three techniques, in Fig. 7. In addition, the deformation of the quarter of the rubber sheet at the maximum load configuration obtained from Technique 3 is illustrated in Fig. 8.

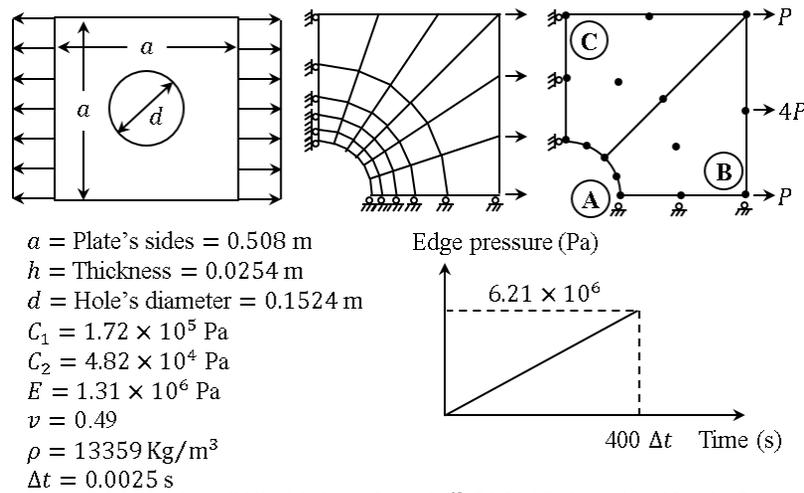

Fig. 6. Left: Geometry and material properties; Middle: Mesh employed in [19]; Right: Mesh employed herein; Bottom: Loading rates. Δt is the average time step for Techniques 1, 2 and 3.

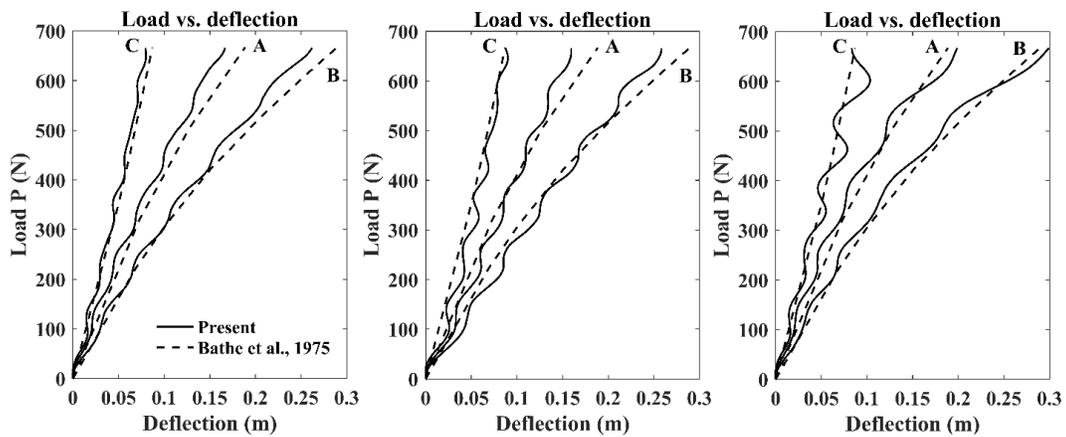

Fig. 7. Load vs. displacement curves at locations A, B and C, obtained from Techniques 1 to 3 from Left to right.

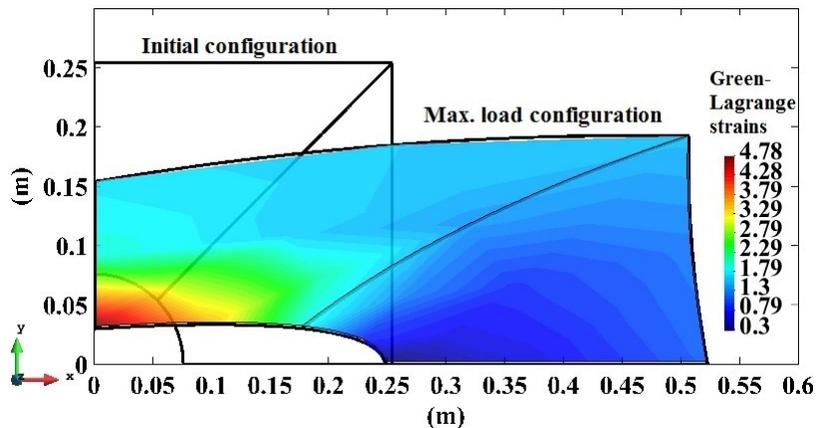

Fig. 8. Illustration of the initial configuration, maximum load (deformed) configuration and the distribution of the Green-Lagrange strains at the max load configuration, obtained from Technique 3.

### 3.3. *Experiment 3: Nonlinear anisotropic incompressible hyperelastic, planar equibiaxial testing, verification of Technique 1*

In [22], square samples of Supple PeriGuard cardiovascular patch material, whose sides were aligned with or orthogonal to the direction of fiber reinforcements, were considered as one of the many potential replacement aortic valve cusp materials. The specimen were subjected to biaxial in-plane tension, from which the average experimental tension vs. displacement data were registered. These authors used the experimental data to determine the material constant for a 2D form of Guccione's material model: $W = \frac{C_1}{2}[\exp(Q) - 1]$, where $Q = C_2 {E_{11}^l}^2 + C_3 \left({E_{22}^l}^2 + \frac{1}{4}\left(\frac{1}{\Delta} - 1\right)^2\right) + 2C_4 {E_{12}^l}^2$, with $\Delta = (2E_{11}^l + 1)(2E_{22}^l + 1)$. We simulated the equibiaxial tension analysis of one quarter of a sample using a $1 \times 1$ mesh of the present thick CB shell FE (Fig. 9), using the material constants and material model reported in [22]. In addition, we used a 3D form of Guccione's material model, where $Q = C_2 {E_{11}^l}^2 + C_3 \left({E_{33}^l}^2 + {E_{22}^l}^2 + {E_{23}^l}^2 + {E_{32}^l}^2\right) + C_4 \left({E_{13}^l}^2 + {E_{31}^l}^2 + {E_{12}^l}^2 + {E_{21}^l}^2\right)$, with the same constants, to verify the accuracy of the present thick CB shell FE in modeling equibiaxial planar stretching of the thin plate (because $h \ll 0.2a$), where transverse shears vanish. The experimental, as well as the present 2D and 3D second Piola-Kirchhoff ($2^{nd}$ PK) membrane tension vs. strain curves are illustrated in Fig. 10. The $2^{nd}$ PK membrane tensions were obtained by multiplying the corresponding components of the second Piola-Kirchhoff stress tensor by the initial thickness of the sample.

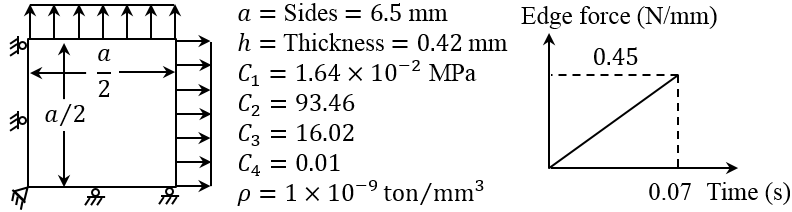

Fig. 9. Geometry, material properties, and loading condition for Experiment 3.

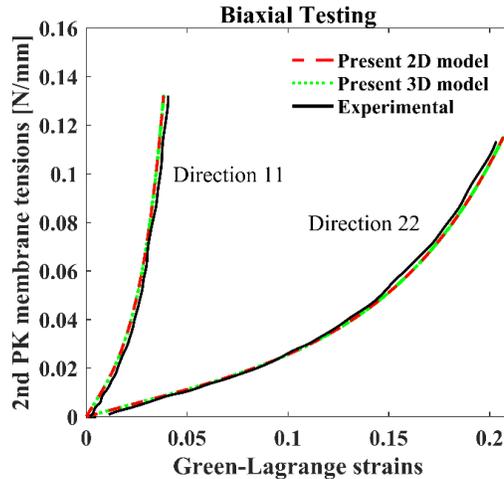

Fig. 10. Force-displacement curves for the equibiaxial in-plane testing, obtained from the 2D and 3D Guccione's material model. The experimental data are taken from Labrosse et al., 2016. Directions 11 and 22 refer to the fiber and cross-fiber directions, respectively.

### 3.4. *Large strain analysis of anisotropic hyperelastic incompressible materials: insensitivity to initially curved geometry and large 3D deformations; verification of Technique 1 and of the fiber length update algorithm*

#### 3.4.1. *Experiment 4: Anisotropic incompressible hyperelastic, cylindrical structure: human thoracic aorta*

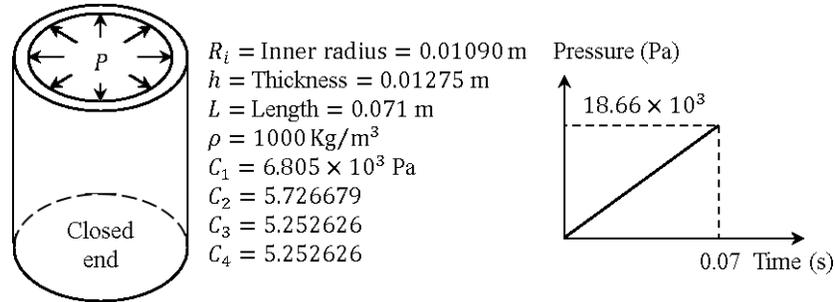

Fig. 11. Left: Human thoracic aorta geometry and material properties; Right: Loading rate.

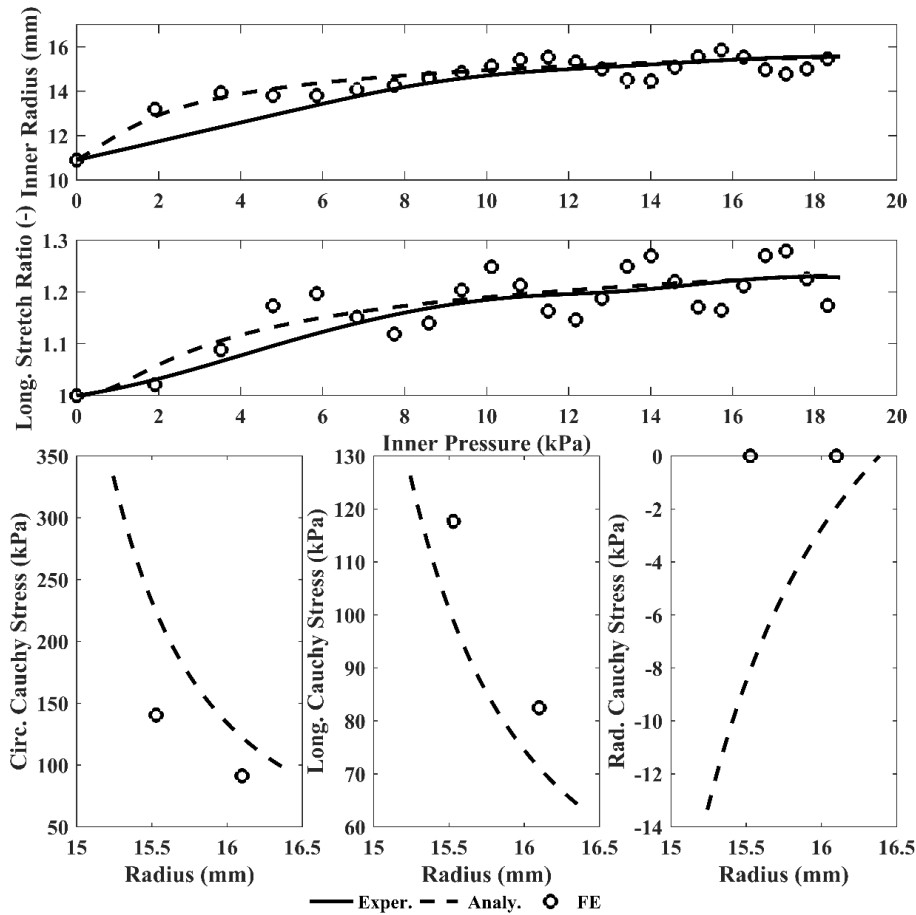

Fig. 12. Comparison between the experimental, analytical and finite element data for human thoracic aorta under pressurization with closed-end and free extension conditions: inner radius vs. pressure (top) and longitudinal stretch ratio vs. pressure (middle). Comparison between the analytical and finite element stress results across the aorta wall in the circumferential (bottom-left), longitudinal (bottom-middle), and radial (bottom-right) directions, measured at 13.33 kPa and close to the open end of the aorta.

In reference [23], segments of fresh human thoracic aortas from eight male sexagenarians were pressurized under closed-end and free extension conditions. Then, material constants for different 3D hyperelastic anisotropic constitutive models were determined from the experimental data. In addition, they modeled the aortic segments in commercial FE code LS-Dyna (LSTC, Livermore, CA, USA), using 6 brick elements in the thickness direction, 20 in the circumferential direction, and about 10 elements in the longitudinal direction. In this study, we used Guccione's 3D hyperelastic anisotropic material model of the form:

$$W = \frac{C_1}{2}\left[\exp\left(C_2 E_{\theta\theta}^{l\,2} + C_3\left(E_{zz}^{l\,2} + E_{rr}^{l\,2} + E_{rz}^{l\,2} + E_{zr}^{l\,2}\right) + C_4\left(E_{\theta z}^{l\,2} + E_{z\theta}^{l\,2} + E_{\theta r}^{l\,2} + E_{r\theta}^{l\,2}\right)\right) - 1\right], \quad (11)$$

where the material constants (Fig. 11, left) were taken from [23]. We used Technique 1 (Section 2.4.1) to develop the constitutive relations, and enforce the material incompressibility and the plane stress condition (in the radial direction presented by index $rr$). In the present FE model, we pressurized one quarter of the cylindrical structure from 0 to 18.66 kPa (physiological range) in 0.07 seconds (Fig. 11, right), using 2 of the present CB shell FEs in the circumferential directions, and 8 elements in the longitudinal direction. The present pressurization rate corresponds to the physiological pressurization rate in early systole, that is 80 mmHg in 0.04 seconds [24]. Comparison between the experimental, analytical and FE data is presented in Fig. 12. To verify that the stress discontinuities (stress jumps) were small in the structure using the present 2×8 mesh, we plotted pressure bands over the deformed geometry at the maximum load configuration. Pressure band (12) is a more conclusive way of presenting stress discontinuities when 3D stresses and 3D deformations are present [14].

$$\text{Pressure band} = \frac{-(\sigma_{xx} + \sigma_{yy} + \sigma_{zz})}{3} \quad (12)$$

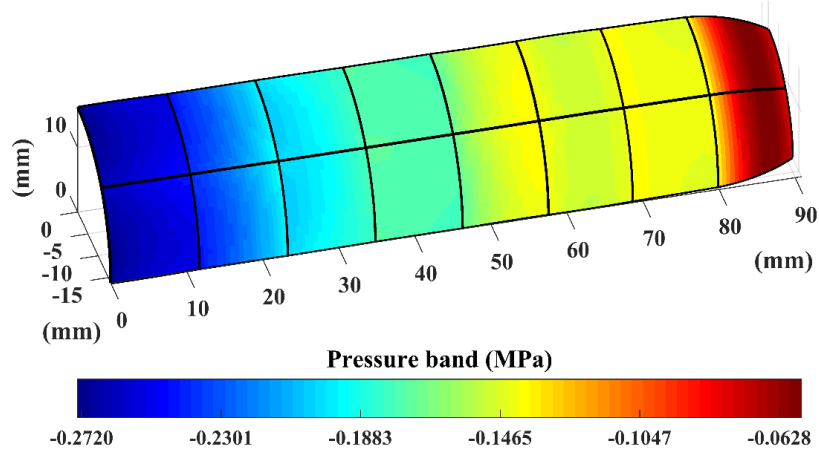

Fig. 13. Deformation and distribution of pressure band (12) due to pressurization under closed-end and free extension conditions at the maximum load configuration.

3.4.2. *Experiment 5: Anisotropic incompressible hyperelastic, cylindrical structure: dog carotid artery*

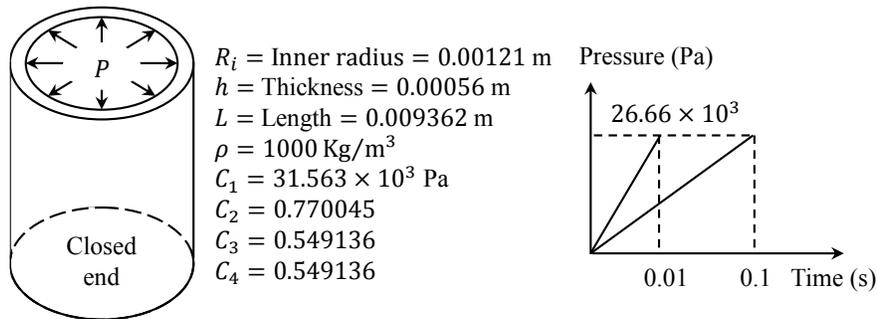

Fig. 14. Left: Dog carotid artery geometry and material properties; Right: Loading rates.

In this experiment, we studied the pressurization of a segment of a closed-end dog carotid artery. The model consisted of an axisymmetric, thick-walled, closed-end cylinder (Fig. 14, left) with homogeneous hyperelastic anisotropic Guccione's material model presented by (11), where we obtained the material constants from the experimental data published in [25] using an approach similar to that presented in [23]. The structure was pressurized from $0 - 26.66$ kPa at two different loading rates (Fig. 14, right). The slower loading rate (26.66 kPa in 0.1 second) was consistent with the physiological pressurization rate before the cardiac cycle starts in early systole, and the faster loading rate (26.66 kPa in 0.01 second) was chosen to study the effect of loading rate on the amplitude of the oscillations in the solution. Comparison between the experimental, analytical and finite element data obtained from the two loading rates is presented in Fig. 15. To verify that the stress discontinuities were small everywhere in the structure, we plotted pressure bands over the deformed geometry at the maximum load configuration obtained from the faster (Fig. 16) and the slower (Fig. 17) loading rates.

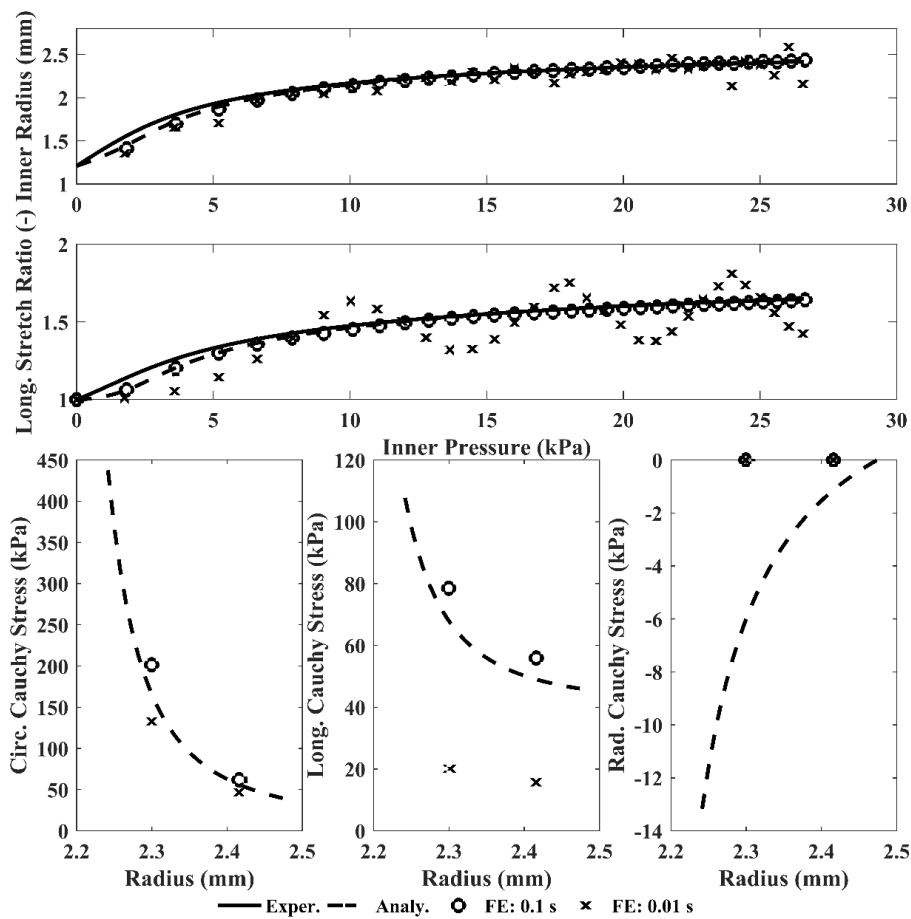

Fig. 15. Comparison between the experimental, analytical and finite element data for dog carotid artery under pressurization with closed-end and free extension conditions: inner radius vs. pressure (top) and longitudinal stretch ratio vs. pressure (middle). Comparison between the analytical and finite element stress results across the aortic wall in the circumferential (bottom-left), longitudinal (bottom-middle), and radial (bottom-right) directions, measured at 13.33 kPa and close to the open end of the artery.

## 4.  Results and Discussion

### 4.1.  *Experiments*

In Experiment 1, 15 of the 9-noded CB shell FE presented in Part 1 were used to model pure bending deformation of an initially flat cantilevered element. As expected, Hughes' fiber coordinate system induced twisting in the cantilevered beam, where pure tip bending moment was applied (Fig. 3). The good accuracy between the present FE model and the analytical solution for all range of bending deformations ($0° − 360°$), illustrated in Fig. 4, is an indication of the reliability of the present algorithm for updating the fiber coordinate system, as well as the insensitivity of the present CB shell FE to shear and membrane locking even though very large bending deformations (360°) are considered.

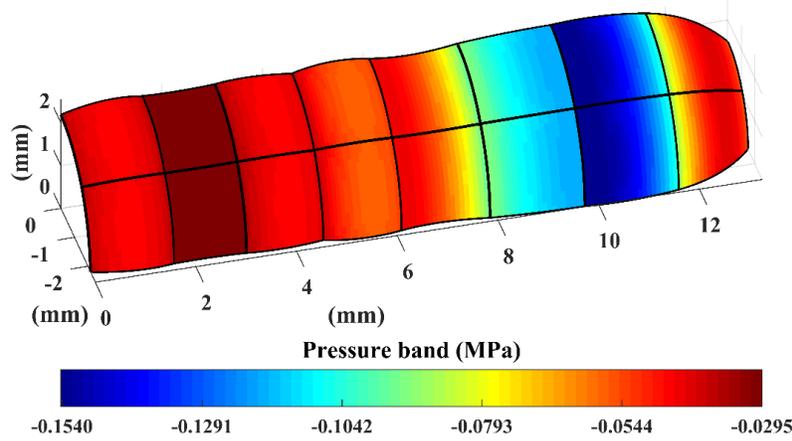

Fig. 16. Deformation and distribution of pressure band due to pressurization under closed-end and free extension conditions at the maximum load configuration obtained from the faster loading rate (26.66 kPa in 0.01 second).

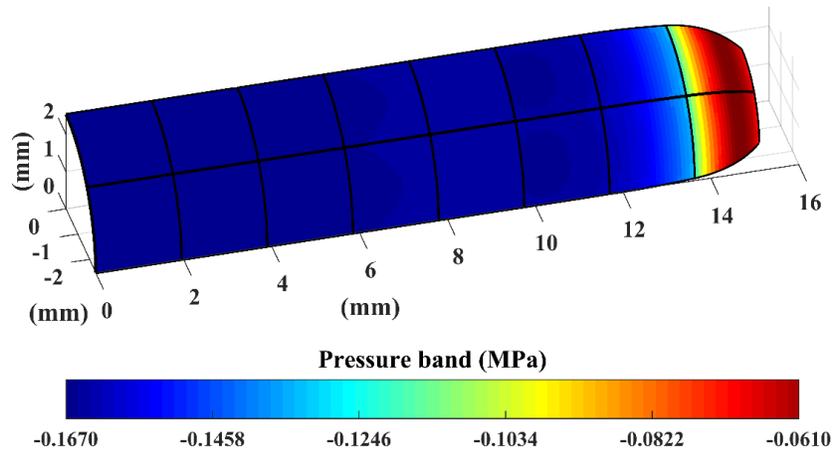

Fig. 17. Deformation and distribution of pressure band due to pressurization under closed-end and free extension conditions at the maximum load configuration obtained from the slower loading rate (26.66 kPa in 0.1 second).

To verify the correctness of the three techniques presented in Section 2.4 for large strain analysis, together with the accuracy and efficiency of our 9-noded CB shell finite element when subjected to large in-plane distortions, as well as the use of an initially irregular mesh and geometry, we performed Experiment 2. In this test, we used the Mooney-Rivlin (hyperelastic, invariant based) material model for Technique 1 (Section 2.4.1), and a constant isotropic constitutive tensor for Technique 2 (Section 2.4.2) and Technique 3 (Section 2.4.3). Basically, Techniques 2 and 3 are designed to work with rubber-like materials that can be described by Young's modulus and Poisson's ratio (i.e. a constant constitutive tensor). To enable comparison between the present explicit dynamic analysis and the static analysis of [19], we ramp-loaded the rubber sheet at an arbitrary 400 increments to the final load (equivalent to $6.21 \times 10^6$ Pa/s) (Fig. 6, bottom-right). The total UL load deformation responses at points A, B, and C obtained from Technique 1 using the Mooney-Rivlin model (Fig. 7, left), is in good agreement with those obtained from Technique 2 using the constant isotropic constitutive tensor (Fig. 7, middle). In addition, the results obtained from these techniques follow the deformation responses at the specified static load increments described in [19]. Carrying out the incremental UL formulation via Technique 3 (using the same constant isotropic constitutive tensor as in Technique 2), the load vs. displacement curves obtained for points A, B, and C were in excellent agreement with those of [19] (Fig. 7, right). Performance-wise, due to the bypass of the transformation of the fourth-order constitutive tensor (done through eight nested for-loops), Techniques 2 and 3 are computationally less expensive than Technique 1. The accuracy of the responses obtained proves all the three techniques valid and accurate in modelling large strains. At the maximum load configuration (illustrated in Fig. 8), Green-Lagrange strains of up to 4.8 are measured, which also numerically proves that the strains analysed are indeed very large, even though a linear constitutive tensor was used. This accuracy in response was achieved using only two of the present CB shell FE, as opposed to the thirty plates used in [19]. Thus, the present CB shell FE is proven insensitive to the adoption of initially irregular elements and large distortions (Fig. 8). Furthermore, the volume of the rubber sheet changed by 0.5% only with such a large in-plane stretching, thereby validating the nodal fiber length adaptation algorithm presented in Part 1 (Section 2.6).

The accuracy of Techniques 1 to 3 in a flat and irregular geometry of a rubber-like material undergoing large in-plane strain and large in-plane distortions was verified by Experiment 2. Moving closer to the biomedical applications of the present work, we further validated the present CB shell FE and Technique 1 (Section 2.4.1) for planar equibiaxial testing of a square sample of Supple PeriGuard (a potential replacement material for aortic valve cusp) using a $1 \times 1$ mesh, and the 2D and 3D Guccione's anisotropic hyperelastic material models (Experiment 3). Fig. 10 exhibited nonlinear and anisotropic tension vs. strain behavior. The 11 direction corresponds to the orientation of the collagen fibers, and the 22 direction is transverse to the collagen fibers. As illustrated, the curves exhibited large deformations under small forces, and then underwent proportionally smaller deformations under larger forces. Due to the alignment of collagen fibers, the 11 direction was stiffer than the 22 direction. The FE results obtained from the 2D and 3D Guccione's models were almost identical to the experimental data obtained in [22]. This proved the accuracy of Technique 1 for 2D and 3D material models, and the robustness of the present thick CB shell FE (accounting for transverse shears) in planar equibiaxial straining of thin structures, in which transverse shears vanish.

Moving forward with the biomedical related applications, we considered modeling pressurization of human thoracic aorta (Experiment 4) and dog carotid artery (Experiment 5) with closed-end and free extension conditions. The aim of these experiments was to further verify Technique 1, the procedure of direct application of incompressibility and plane stress conditions, and the fiber length update algorithm for large 3D strains (Part 1, Section 2.6) in thin (Experiment 4) and thick (Experiment 5) shell structures, when initially curved geometries (uncoincident lamina and global coordinate system) are considered, and the structure undergoes large 3D deformations (causing large discrepancies between lamina and fiber coordinate systems). These experiments confirmed that the presence of transverse shears in the formulation of the present thick CB shell FE does not cause artificial stiffening in large 3D deformations of a thin structures, where transverse shears are physically not present.

We used Guccione's 3D hyperelastic anisotropic material model (11) for both experiments. Note that the Guccione's anisotropic hyperelastic material model was selected because it is able to handle full 3D strains, and because we could evaluate the corresponding material constants. It is worth emphasizing the very important point that 2D or membrane type material models include in-plane strains only, and cannot describe 3D material properties. In particular, features such as residual stresses, through-thickness stress distributions, and torsional deformations cannot be captured with a 2D material model. Comparisons between the experimental, analytical and finite element results conducted for Experiments 4 and 5 are, respectively, shown in Fig. 12 and Fig. 15 for the inner radius (top) and the longitudinal stretch ratio (middle) upon pressurization. For practical purposes, the curves exhibit excellent agreements between each other, over the specified pressure ranges. Considering Experiment 4, the FE results presented in [23] were obtained from LS-Dyna using 6 brick elements in the thickness direction, 20 in the circumferential direction, and 10 elements in the longitudinal direction (i.e. $6 \times 20 \times 10$ mesh). This, in comparison with the present $2 \times 8$ mesh, proves the present CB shell FE more efficient than the brick elements. In addition, the critical time step evaluated by LS-Dyna was in the range of $6 \times 10^{-8}$ to $10 \times 10^{-8}$ seconds, whereas the present average critical delta time evaluated from the eigen value approach presented in Part 1 (Section 2.9) was $1 \times 10^{-5}$ seconds, which is about $8 \times 10^3$ times larger than that of LS-Dyna.

Considering that in the formulation of our CB shell FE, the reference surface of the shell is taken to be the mid-surface, the analysis is done on the mid-surface (middle radius) of the structure and the values of the internal radius are obtained from subtracting half of the thickness (nodal fiber lengths) from the middle radius. Thus, the accuracy of the present FE inner radius vs. pressure curves (Fig. 12 and Fig. 15, top), as well as the negligible change in the volume of the structure (−0.8% and −0.3% for Experiments 4 and 5, respectively) confirmed the accuracy of the fiber length update algorithm for large 3D strains and large 3D distortions of curved geometries. The more pronounced oscillatory response in the longitudinal direction (Fig. 12 and Fig. 15, top and middle) was consistent with the anisotropic behaviour the aorta, where the longitudinal direction is softer than the circumferential direction. Comparison between the analytical and the present FE stress distributions across the wall (pressurized to 13.33 kPa) in the three circumferential, longitudinal, and radial directions is shown in Fig. 12 and Fig. 15 (Experiments 4 and 5, respectively), bottom- left to right, in that order. The analytical curves illustrate that the stresses were the highest towards the inside of the aortic wall for the circumferential and longitudinal directions, as expected in the absence of residual stresses in the model. As illustrated in Fig. 15 (top and middle), the slower loading rate decreased the amplitude of the oscillations, resulting in a more favorable comparison between the analytical and the present FE stress distributions (Fig. 15 bottom-middle). Finally, the analytical radial stresses were the highest towards the inside of the wall (Fig. 10 and Fig. 13 bottom-right), where the stress magnitude was equal to that of the applied pressure (13.33 kPa), and equal to zero on the outside surface, where no external pressure was applied. Consistent with the plane stresses assumption associated with the Mindlin-Reissner shell theory (thus, the present CB shell FE), the radial stresses (normal to the shell surface) obtained from the present CB shell FE remained zero across the element thickness. Comparison between the analytical circumferential and longitudinal (in-plane) stresses, and the analytical radial (normal) stresses validated that, because normal stresses are much smaller than the in-plane stresses, they can be neglected in structural theories.

      Comparison between Fig. 16 and Fig. 17 suggested that not only decreasing the loading rate decreased the amplitudes of the response oscillations (Fig. 15), but also produced a smoother deformed shape (i.e. removed the wrinkling along the vessel segment), and reduced stress jumps between the elements (i.e. improved stress convergence). Overall, the stress jumps, although small across each element, appeared to be the most pronounced in the elements closest to the closed end. This was because those elements underwent the largest 3D deformations and distortions. Should it be required, a smoother pressure band could be achieved either through mesh refinement (as verified in Part 1), or by decreasing the loading rate (as verified trough Experiment 5).

      In conclusion, the results obtained from Experiments 4 and 5 validate Technique 1, the procedure of direct application of incompressibility and plane stress condition (Section 2.4.1), as well as the fiber length update algorithm for large membrane strains in thin and thick shell structures (Part 1, Section 2.6). In addition, the present CB shell FE was shown to be more efficient and accurate compared to the brick elements (e.g. as used in LS-Dyna [3]), since much fewer elements and much larger critical time step were needed to accurately model very large 3D deformations and very large 3D strains of thin and thick curved structures.

### 4.2. *Contributions*

In Part 1, the present thick CB shell FE, based on the second order shear deformation theory (Mindlin-Reissner) was shown to be a serious competitor to higher order shell elements, requiring fewer (5) degrees of freedom per node, or iterations to achieve displacement and stress convergence. The CB shell FE is also naturally shear, membrane and volumetric locking insensitive, and allows for relatively large critical time step ($\Delta t$). In addition, the present fiber coordinate system algorithm is insensitive to large (360°) bending deformations (Experiment 1).

      The geometric and kinematic description of the present CB shell FE combined with the present kinetic descriptions and the three techniques to evaluate the hyperplastic constitutive relations in the lamina coordinate system, enabled the accurate modeling of very large 3D deformations and very large 3D strains with a coarser mesh of initially irregular elements and geometries. Thus, the present work is applicable to rubber-like materials and soft biological tissues (Experiments 2 to 5).

### 4.3. *Limitations*

The present 9-noded thick CB shell FE was proven efficient and accurate, not only in modeling linearly elastic materials (engineering applications presented in Part 1), but also in modeling anisotropic incompressible hyperelastic materials such as rubber-like materials and soft biological tissues (presented herein). To make a full use of this element in surgical simulations (e.g. heart valve repair modeling), the formulations still need to be expanded to integrate contact constraints in the finite element analysis. Contact problems are categorized as boundary nonlinearities, because both contact boundaries and contact stresses are unknown, and there is an abrupt change in contact forces. Many contact algorithms exist in the literature [6,14,26]. In addition, for the study of arteries, it will be important to incorporate residual stresses. The methods to do so are similar whether for brick or shell elements. The issue of inhomogeneous layers can be handled by overlaying multiple CB shell elements and assigning them different material properties.

## 5. Conclusion

Geometrically and materially nonlinear FE analysis requires a CB shell FE that is accurate, reliable and versatile. The accuracy and efficiency of the proposed CB shell FE for large deformations of linear elastic materials was tested in Part 1. In the present study, we further verified these qualities in large 3D deformations and large 3D strains of hyperelastic anisotropic incompressible materials (rubber-like materials and soft biological tissues). The CB shell

FE proved to be a highly accurate, reliable, and efficient tool for modeling, analyzing, and predicting the mechanical behavior of soft biological tissues in physiological loading conditions. Overall, the present work is a promising step toward the robust simulation of soft biological tissues such as heart valves and blood vessels.

**Acknowledgements**

This work was supported by the Natural Sciences and Engineering Research Council of Canada for Discovery Grant 312065-2012 (M.R.L).